\documentclass[12pt]{amsart}
\theoremstyle{plain}
\newtheorem{theorem}{Theorem}[section]
\newtheorem{corollary}[theorem]{Corollary}

\newtheorem{proposition}[theorem]{Proposition}
\newtheorem{lemma}[theorem]{Lemma}

\theoremstyle{definition}
\newtheorem{definition}[theorem]{Definition}

\theoremstyle{remark}

\renewcommand{\a}{\alpha}
\renewcommand{\b}{\beta}
\renewcommand{\d}{\delta}

\newcommand{\e}{\varepsilon}
\newcommand{\g}{\gamma}
\newcommand{\G}{\Gamma}

\newcommand{\var}{\varphi}
\newcommand{\s}{\sigma}

\newcommand{\bz}{\mathbb Z}
\newcommand{\bzn}{{\mathbb Z}^n}
\newcommand{\br}{\mathbb R}
\newcommand{\brn}{{\mathbb R}^n}
\newcommand{\bt}{\mathbb T}
\newcommand{\btn}{{\mathbb T}^n}
\newcommand{\ca}{\mathcal A}
\newcommand{\cb}{\mathcal B}
\newcommand{\mcr}{\mathcal R}
\newcommand{\cv}{\mathcal V}

\newcommand{\x}{\times}
\newcommand{\<}{\langle}
\renewcommand{\>}{\rangle}

\begin{document}

\title{Projective multi-resolution analyses for $L^2(\br^2)$}
\author[J.A.~Packer]{Judith A. Packer}
\address{Department of Mathematics, National University
of Singapore, 2 Science Drive 2, Singapore 117543,
Republic of Singapore}
\address{Department of Mathematics, University of
Colorado, CB 395, Boulder, Colorado, 80309-0395, U.S.A.}
\email{packer@euclid.colorado.edu}
\author[M.~A.~Rieffel]{Marc A.~Rieffel}
\address{Department of Mathematics, University of
California at Berkeley, Berkeley, CA 94720, U.S.A.}
\email{rieffel@math.berkeley.edu}
\keywords{wavelets, multi-resolution, tight frames,
projective modules, module frames,
Hilbert C*-module, $K$-theory}
\subjclass{Primary 46L99; Secondary 42C15, 46H25, 47A05}
\thanks  {The first author was supported  in part by a
National University of Singapore Research Grant R-146-000-025-112.
The research of the second author was supported in part
by National Science Foundation grants DMS99-70509
and DMS-0200591.}

\begin{abstract}
   We define the notion
of  ``projective" multiresolution analyses, for which,
by definition, the initial space corresponds to a
finitely generated projective module
over the algebra $C(\btn)$ of continuous
complex-valued functions on an $n$-torus. The case of ordinary
multi-wavelets is that in which the projective module is
actually free. We discuss the properties of projective
multiresolution analyses, including the frames which they
provide for $L^2(\brn)$. Then we show how to construct
examples for the case of any
diagonal $2 \times 2$ dilation matrix with integer entries,
with initial module
specified to be any fixed finitely generated
projective $C(\mathbb T^2)$-module. We compute the
isomorphism classes of the corresponding wavelet modules.

\end{abstract}
\maketitle

In classical wavelet theory one uses
multi-resolution analyses to construct
(multi-) wavelets and their corresponding orthonormal bases
or frames for $L^2(\brn)$. In almost all applications the
scaling functions and wavelets have continuous Fourier
transforms. This continuity is a significant and interesting
condition. If one requires just a bit more, then one finds
that in the frequency domain one is dealing with what are
called projective modules over $C(\btn)$ (or equivalently,
with the spaces of continuous cross-sections of complex
vector bundles over $\btn$). The case of a single scaling
function or wavelet corresponds to the free module of rank
1, whereas the case of several orthogonal scaling functions
or wavelets corresponds to free modules of higher rank. This
leads one to ask whether wavelet theory carries over to the
case of general projective modules over $C(\btn)$. It is the
purpose of this paper to show that the answer is affirmative.

For a given dilation matrix $A$ we define a projective
multiresolution analysis for $L^2(\brn)$ to be an increasing
sequence $\{V_j\}$ of subspaces having the usual properties,
with the one exception that instead of $V_0$ being the linear
span of the integer translates of one or more scaling functions,
we only require that $V_0$ be the (projective) module of
inverse Fourier transforms of continuous
cross-sections of a complex vector bundle over $\btn$.
The precise definition
is given in Definition \ref{defpr}. We show in
Section 3 that projective
multi-resolution analyses give in a natural way wavelets which
determine normalized tight frames for $L^2(\brn)$. This seems
to be a somewhat new way of constructing tight frames
for $L^2(\brn)$. We note that the term ``projective multiresolution analysis" 
was first used by G. Zimmermann his thesis \cite{Z1}; the definition given by 
him is different from the one used here, and we comment on the difference in Section 3.

It is not immediately evident that non-trivial
projective multiresolution
analyses exist. We show how to construct them for any
diagonal $2 \times 2$ dilation matrix and the continuous
cross-section space of any complex vector bundle over $\bt^2$.
This fleshes out ideas which the second author presented in
a special session talk at the annual meeting of the AMS
in 1997 \cite{Rie2}. We also determine the isomorphism
class of the corresponding wavelet modules. We find that
in some cases the wavelet modules are free modules, while
in other cases they are not free.

As hinted above, we
find it very convenient to take a fairly algebraic
approach to our topic. The function spaces determined by scaling
functions and wavelets are viewed as modules over the convolution
algebras of the corresponding discrete subgroups of $\brn$.
In the frequency domain these convolution algebras are
(when suitably completed) isomorphic to $C(\btn)$, which is a commutative
$C^*$-algebra. As seen already in \cite{pacrief}, it is very
useful to define on modules over $C(\btn)$ inner products with
values in this algebra. Modules equipped with algebra-valued
inner products have been used very profitably for three decades
in related areas of non-commutative harmonic
analysis \cite{Rie0} \cite {La}, and
are often called ``Hilbert $C^*$-modules''. Recently there has
been much interest in frames for Hilbert $C^*$-modules.
Besides \cite{pacrief} see \cite{LaFr98} \cite{LaFr00} for
discussions and references. Such frames turn out to be very
useful in our context, leading naturally to ordinary
tight frames for $L^2(\brn)$.

Our paper is organized as follows. In Section 1 we
discuss the setting for our topic,
and describe in detail the standard $C(\btn)$-module
$\Xi$ where we will perform most of our constructions.
Roughly speaking, $\Xi$ is the closure of $C_c(\brn)$
in $L^2(\brn)$ for a norm stronger than the Hilbert-space
norm. In Section 2 we recall the definition and basic
properties of projective modules, and of the frames which
are associated with them. In Section 3 we define the notion of projective
multiresolution analyses, and obtain many of their general
properties. In particular, we discuss the module frames
and the tight frames for $L^2(\brn)$ associated with them.
In Section 4
we briefly review specific constructions of certain finitely
generated projective $C(\mathbb T^2)$ modules along the
lines discusses
in \cite{Rie1}.  They provide model modules for each isomorphism
class of projective modules over $C(\bt^2)$.
We mention some of their properties, and
their relationship to the $K$-theory of $\bt^2.$

In Section 5 we restrict attention to diagonal $2 \times 2$
dilation matrices. We present there one of our main technical results, namely,
the construction of an appropriate substitute for a
scaling function in $\Xi$ corresponding to a desired initial
projective module.
We give a sufficient condition which
guarantees that an initial $C(\bt^2)$-submodule of $\Xi$
generated by a
single element of $\Xi$ will gives rise to a projective
multi-resolution analysis corresponding to a given dilation operator.
We use this to show how to construct, for a given arbitrary
projective $C(\bt^2)$-module and an arbitrary diagonal
dilation matrix, a projective multi-resolution analysis
with its initial module $V_0$
isomorphic to the given projective module.
In Section 6 we identify the structure and isomorphism classes of the
wavelet modules corresponding to the projective
multiresolution analyses constructed
in the Section 5.

We should note that for diagonal matrices, some our
constructions extend to higher dimensions, and we intend
to publish these results in a later paper.  For matrices
which are not similar via an element of $GL(2,\mathbb Z)$
to a diagonal matrix, the situation seems more difficult;
however we note that we have been able to constructed a projective
multi-resolution analysis for the matrix
$A\;=\;\left(\begin{array}{rr}
1&-1\\
1&1
\end{array}\right)$ whose initial module is not free.  The methods used in our
construction at this point seem somewhat ad-hoc; but one
intriguing result from the calculation in this case is that
the higher dimensional modules $W_i$ are all free
$C(\mathbb T^2)$-modules of dimension $2^{i+1}-2^{i}$, which we
conjecture corresponds to the determinant of $A$ being
positive, as in the diagonal case.

The first author would like to thank Professors Lawrence
Baggett and Dana Williams for many useful discussions on
the topics discussed in this paper, and for their great
hospitality towards her and her family during her visit
to the University of Colorado at Boulder and Dartmouth
College, respectively, during her sabbatical year.

\section{The setting}

In this section we describe the setting within which our
construction of projective multiresolution analyses will
take place. Much of the material in the first part of this section is
well-known, but it serves to establish our notation.
We first sketch our setting in the ``time'' domain,
but then we indicate more precisely what the setting
looks like in the ``frequency'' domain, where, as
usual in this subject, matters are more transparent, and
where we will work for the rest of the paper.

We will denote the algebra of bounded operators on $L^2(\brn)$
by ${\mathcal L}(L^2(\brn))$.
For each $s \in \brn$ we define the translation operator
$T_s$ on $L^2(\brn)$ by $(T_s\xi)(t) = \xi(t - s)$.
We denote by $\ca$ the operator-norm closed subalgebra
of the ${\mathcal L}(L^2(\brn))$ generated by
all the translation operators $T_\g$ for $\g \in \bzn$.
Equivalently, $\ca$ is generated by all the convolution
operators $T_f$ for $f \in \ell^1(\bzn)$. We will see
below that the elements of $\ca$ will in our setting
act somewhat like scalars, and that as is often done
in related situations \cite{La} \cite{ pacrief} it is then convenient to
use right-module notation, while putting other operators
on the left. Thus for $f \in L^1(\bzn)$
and $\xi \in L^2(\brn)$ we have
\[
(T_f\xi)(t) = (\xi * f)(t) =
\sum_{\g \in \bzn} f(\g)(T_\g\xi)(t) =
\sum_{\g \in \bzn} \xi(t - \g)f(\g).
\]
We will denote by $\| f \|_\ca$ the operator norm of
$T_f$. The natural involution coming from taking
adjoints of operators is defined by $(f^*)(\g) =
{\bar f}(-\g)$.

Throughout we let $A$ be an
$n\times n$ dilation matrix, that is, an invertible
matrix with
integer entries, all of
whose eigenvalues have modulus strictly
greater than $1$. Note that $A\bzn \subset \bzn$.
Throughout we let $\d = |\det(A)|^{-1/2}$. We let $D$
denote the unitary dilation operator on $L^2(\brn)$
defined by
\[
(D\xi)(t) = \d^{-1} \xi(At).
\]
We let $\a$ denote the automorphism of
${\mathcal L}(L^2(\brn))$ consisting of conjugating
by $D$, that is, $\a(T) = DTD^{-1}$ for
$T \in {\mathcal L}(L^2(\brn))$. Now $\a$ does not
carry $\ca$ into itself, but rather onto the algebra
of convolution operators determined by $\ell^1(A^{-1}\bzn)$.
To be specific,
\[
(\a(f))(\g) = f(A\g)
\]
for $f \in \ell^1(\bzn)$ and $\g \in A^{-1}\bzn \subseteq \brn$.
More generally, for any
$j \in \bz$ let $\ca_j$
denote the algebra of operators corresponding to
$\ell^1(A^{-j}\bzn)$. Then $\ca_j \subset \ca_{j+1}$
and $\a(\ca_j) = \ca_{j+1}$, while $\ca_j = \a^j(\ca)$,
for each $j$. These algebras are all isomorphic, but
we must be careful to distinguish between them because they
are quite different sets of operators on $L^2(\brn)$.

As mentioned in the introduction, in order for our
topic to have content it is essential that we work
with functions whose Fourier transforms are continuous.
For this purpose it is natural to introduce an appropriate
dense subspace of $L^2(\brn)$, on which is defined an
$\ca$-valued ``inner product''. This inner product
is defined by convolution
followed by ``down-sampling'', namely
\[\<\xi, \eta\>_\ca(\g) = (\xi^* * \eta)(\g)
= \int_{\brn} {\bar \xi}(t)\eta(t+\g) \ dt.
\]
Notice that this inner product is very natural within
the context of wavelet theory (as seen also in \cite{pacrief}),
since if
$\var_1, \dots , \var_m$ is a possible family of
multi-scaling functions, or multi-wavelets, then the
condition that all their translates by elements of
$\bzn$ form an  orthonormal set is just the condition that
\[
\<\var_j, \var_k\>_\ca = \d_{jk}1_\ca
\]
for all $j, k$, where $1_\ca$ is the ``delta-function''
at $0$ in $\ell^1(\bzn)$. When we pass to the ``frequency
domain'' it will be much clearer what the domain of this
inner product should be. But the Schwartz functions
will be dense in the domain, and so the reader can
temporarily assume that the domain consists of the
Schwartz functions.
It is for this inner product that
the elements of $\ca$ act much like scalars, in the sense that
\[
\<\xi, \eta*f\>_\ca = \<\xi, \eta\>_\ca *f.
\]
In terms of this inner
product we can define an ordinary norm in the same way as
one does for Hilbert spaces, namely
\[
\|\xi\|_\ca = (\|\<\xi, \xi\>_\ca\|_\ca)^{1/2}.
\]

We can define in a similar way inner products on $\Xi$ with
values in each $\ca_j$. Then we will have
\[
\<D\xi, D\eta\>_{\a(\ca_j)} = \a(\<\xi, \eta\>_{\ca_j}),
\]
and similarly for powers of $D$.
Since $\a$ is an isometry for the operator norm,
it follows that $D$
is isometric, in the sense that $\|D\xi\|_{\a(\ca_j)} = \|\xi\|_{\ca_j}$.
It is also easy to check that
\[
D(\xi*f) = (D\xi)(\a(f)).
\]

We now pass to the frequency domain. We define on $\br$
the function  $e$  by $e(r) = e^{2\pi ir}$,
and we define the Fourier transform, $\mathcal F$, by
\[
{\mathcal F}(\xi)(x) =
   \int_{\brn} \xi(t)\bar e(t\cdot x) \ dt.
\]
For the rest of the paper we will work in the frequency
domain, and so for notational simplicity we will make
the unusual choice of denoting the Fourier transforms of
objects by the same notation as the original objects.
Usually only the variables will indicate that we are
in the frequency domain.

Simple calculations show what happens to the setting
described above under the Fourier transform. We
record here some of the facts, many of them familiar.
The algebra $\ca = \ca_0$
becomes the algebra $\ca = C(\btn)$ of continuous
functions on $\btn = \brn/\bzn$. The action of $\ca$
on $L^2(\brn)$ is now given simply by pointwise
multiplication, where we view the functions in $\ca$
as periodic functions on $\brn$. That is,
\[
(\xi f)(x) = \xi(x)f(x).
\]
The operator norm on $\ca$ is now just the usual supremum
norm $\|\cdot\|_\infty$, and the involution is
complex conjugation. We have
\[
({\mathcal F}D{\mathcal F}^{\ast})(x) = \d\xi((A^t)^{-1}x),
\]
where $A^t$ is the transpose of $A$. Set $B = (A^t)^{-1}$.
Then we see that in the frequency domain we deal with
the operator ${\mathcal F}D{\mathcal F}^{\ast}.$ As we will primarily
be working in the frequency domain, there is no danger of confusion
by denoting this operator by $D$ as well, and $D$ is given by
\[
(D\xi)(x) = \d\xi(Bx).
\]
To see what the
$\ca$-valued inner product will be, we first determine the
Fourier transform of the down-sampling operator, which we
now denote by $S$. Let $\xi$ be a Schwartz function in the
time domain. Define $e_x$ by $e_x(t) = e(t\cdot x)$.
By the Poisson summation formula we obtain for
$x \in \brn$
\begin{eqnarray*}
(S\xi)\hat{}(x) &=& \sum_\g \xi(\g) \bar e(\g \cdot x)
= \sum_p (\xi \bar e_x)\hat{}(p) \\
&=& \sum_p\int (\xi \bar e_x)(t)\bar e(t\cdot p) \, dt
= \sum_p \int \xi(t)\bar e(t\cdot x)\bar e(t\cdot p)\, dt
= \sum_p \hat \xi(x+p).
\end{eqnarray*}
Thus for $\eta$ a Schwartz function in the frequency domain
we have
\[
(S\eta)(x) = \sum_p \eta(x+p) = \sum_p \eta(x-p).
\]
In view of this, we see that
the inner product in the frequency domain
is given by
\[
\<\xi, \eta\>_\ca(x) =
\sum_{p \in \bzn}(\bar \xi \eta)(x-p).
\]
It is clear from this that
\[
\<\xi, \eta f\>_\ca = \<\xi, \eta\>_\ca f
\]
for $f \in \ca$.
We can now see more easily what the domain, $\Xi$,
for this inner product
should be.

Various versions of the space $\Xi$ below have been introduced and studied before,
specifically by G. Zimmermann in Chapter V of his thesis \cite{Z1},  where the notation $L^{2,\infty}(\mathbb Z^n,\mathbb T^n)$ was used, again by Zimmermann in \cite{Z2}, and by J. Benedetto and Zimmermann in \cite{BZ}.  We thank the referee for bringing these works to our attention. We use different notation here, which we specify in the following definition.
  
\begin{definition}
\label{def1a}
We let $\Xi$ denote the set of bounded continuous functions
on $\brn$ for which there is a
constant, $K$, such that $\sum_{p \in \bzn}|\xi(x-p)|^2 \leq K$
for each $x \in \brn$, and furthermore such that the
function defined by this sum is continuous. The norm
is given by
\[\|\xi\|_\ca = \|\<\xi,\xi\>_\ca\|^{1/2}
= \sup_x(\sum_{p \in \bzn}|\xi(x-p)|^2)^{1/2}.
\]
\end{definition}

We remark that when one sees the naturalness of the proofs of
Propositions \ref{prop2i} and \ref{prop2d} one understands
that $\Xi$ with its $\ca$-valued inner product is a quite
comfortable setting within which to develop even ordinary
wavelet theory.

\begin{proposition}
\label{prop1a}
For any $\xi, \eta \in \Xi$ the sum defining $\<\xi,\eta\>_\ca$
converges uniformly on compact subsets of $\brn$, and
we have
\[
\<\xi, \eta\>_{L^2(\brn)} = \int_{\btn}\<\xi,\eta\>_\ca \ \ .
\]
For any $\xi \in \Xi$ we have
\[
\|\xi\|_2 \leq \|\xi\|_\ca \quad \quad and \quad\quad
\|\xi\|_\infty \leq \|\xi\|_\ca.
\]
\end{proposition}
\begin{proof}
Dini's lemma (page
34 of \cite{Ped}) tells us that the sum for $\<\xi, \xi\>_\ca$
will converge
uniformly on compact subsets of $\brn$. The corresponding
fact for $\<\xi, \eta\>_\ca$ then follows by polarization.
   For any $\xi, \eta \in \Xi$ we
have
\begin{equation*}
\<\xi,\eta\>_{L^2(\brn)} = \int_{\brn} \bar \xi(x)\eta(x)\ dx =
\int_{\btn}(\sum_{p \in \bzn}\bar\xi(x-p)\eta(x-p))\ dx
= \int_{\btn} \<\xi,\eta\>_\ca(x) \ dx.
\end{equation*}
The first inequality follows immediately from this.
The second inequality is clear from the definitions.
\end{proof}

We will need the following version of the Cauchy-Schwarz
inequality for our inner product. (See proposition 2.9
of \cite{Rie0}, and \cite{La}.)

\begin{proposition}
\label{prop1cs}
For any $\xi, \eta \in \Xi$ we have
\[
|\<\xi, \eta\>_\ca|^2 \leq \< \xi, \xi\>_\ca \<\eta, \eta\>_\ca
\]
as an inequality between functions in $\ca$.
\end{proposition}

\begin{proof}
For any $\zeta \in \Xi$ and any $x \in \brn$ define
$\zeta^x \in \ell^2(\bzn)$ by $\zeta^x(p) = \zeta(x-p)$.
Then for $\xi, \eta \in \Xi$ we have, upon using the ordinary
Cauchy-Schwarz inequality,
\begin{eqnarray*}
|\<\xi, \eta\>_\ca(x)|^2 &=& |\sum_p \bar \xi^x(p)\eta^x(p)|^2 \\
&\leq& (\sum_p |\xi^x(p)|^2)(\sum_q |\eta^x(q)|^2)
=\< \xi, \xi\>_\ca(x) \<\eta, \eta\>_\ca(x).
\end{eqnarray*}
\end{proof}

Standard arguments of the kind used to treat $\ell^2$ show
that $\Xi$ is a linear space on which the values of
the inner product are in $\ca$. We now introduce common
notation which we will use here and in a number of places later.
We let $C_b(\brn)$ denote the algebra of continuous bounded
functions on $\brn$, we let $C_\infty(\brn)$ denote its
subalgebra of functions which vanish at infinity, and we
let $C_c(\brn)$ denote its subalgebra of functions of
compact support. We let $I^n$ denote the cube $[0,1]^n$. Every
coset of $\bzn$ in $\brn$ meets $I^n$. Thus sums such as that
defining $\<\xi,\eta\>_\ca(x)$ need only be considered for
$x \in I^n$.

\begin{proposition}
\label{prop1p}

The space $\Xi$ is complete for its norm.
Furthermore, $\Xi \subset C_\infty(\brn)$.
\end{proposition}

\begin{proof}
Suppose that $\{\xi_k\}$ is a Cauchy sequence in $\Xi$. Then
from Proposition \ref{prop1a} it follows that this sequence
is Cauchy for $\|\cdot\|_\infty$, and so
converges uniformly to some $\xi \in C_b(\brn)$. For any
$\eta \in \Xi$ and any $x \in \brn$ let $\eta^x$ be as in
the proof or Proposition \ref{prop1cs}.
Because $\{\xi_k\}$ is a Cauchy sequence, so is $\{\xi_k^x\}$
for each $x \in \brn$, and $\{\xi_k^x\}$ will converge in
$\ell^2(\bzn)$, necessarily to $\xi^x$.
Furthermore, from the definition of $\|\cdot\|_\ca$ this convergence
is uniform in $x$. It follows that $\|\xi\|_\ca = \sup_x\|\xi^x\|_2$,
so that $\xi \in \Xi$. A bit of further argument shows that
for any $\e > 0$ there is a $k_0$ such that if $k \geq k_0$ then
$\|\xi^x - \xi_k^x\|_2 < \e$ for all $x\in \brn$,
so that $\{\xi_k\}$ converges to $\xi$ for the norm of $\Xi$.

Finally, let
$\xi \in \Xi$ be given.  For any given $\e > 0$ we can find
a finite subset $S$ of $\bzn$ such that $\sum_{p \notin S}|\xi(x-p)|^2 < \e$
for all $x \in I^n$. Then, in particular, $|\xi(x-p)|^2 < \e$ for
every $p \notin S$ and every $x \in I^n$.
Thus if $|\xi(y)|^2 > \e$ for some $y \in \brn$, we have
$y \in I^n - S$ (the union of translates of $I^n$ by
elements of $S$), which is compact. Consequently,
$\xi \in C_\infty(\brn)$.
\end{proof}

\begin{proposition}
\label{prop1c}
   If $F \in C_b(\brn)$
then the pointwise product $\xi F$ is in
$\Xi$ for each $\xi \in \Xi$, and
$\|\xi F\|_\ca \leq \|\xi\|_\ca \|
F\|_\infty$. Furthermore, $C_c(\brn)$ is dense
in $\Xi$ for its norm.
\end{proposition}

\begin{proof}
For any finite subset $S$ of $\bzn$ we have
\[
\sum_S |(\xi F)(x-p)|^2 \leq \|F\|_\infty \sum_S |\xi(x-p)|^2,
\]
and from this it easily follows that $\xi F \in \Xi$ and
$\|\xi F\|_\ca \leq \|\xi\|_\ca \|
F\|_\infty$.

It is easily seen that $C_c(\brn) \subset \Xi$. We now show that
$C_c(\brn)$ is dense. Given $\xi$ and $\e > 0$, choose
a finite set $S \subset \bzn$ such
that $\sum_{p \notin S} |\xi(x-p)|^2 < \e$ for all $x \in I^n$.
Choose an $F \in C_b(\brn)$ such that $F = 0$ on $I^n - S$
   and
$F(y) = 1$ for all $y$ outside some compact neighborhood of
$I^n - S$.
Then $\xi(1-F) \in C_c(\brn)$ and
$\|\xi - \xi(1-F)\|_\ca = \|\xi F\|_\ca < \e$.
\end{proof}

\begin{proposition}
Let $\xi \in C_b(\brn)$. If there is a constant, $K$,
and an $s > n/2$ such that $|\xi(x)| \leq K\|x\|^{-s}$
for all sufficiently large $x$, then $\xi \in \Xi$.
\end{proposition}

\begin{proof}
For $x \in I^n$ and all sufficiently large $p \in \bzn$
we will have $|\xi(x-p)| \leq 2K\|p\|^{-s}$.
Since $\sum_p \|p\|^{-2s} < \infty$, it follows that
$\sum |\xi(x-p)|^2$ converges
uniformly for $x \in I^n$, to a necessarily continuous
function. By translating by $\bzn$ it follows that the
sum converges uniformly on any compact subset of $\brn$,
to a continuous function.
\end{proof}

\begin{corollary}
\label{cor1s}
The space $\Xi$ contains all the Schwartz functions on $\brn$.
\end{corollary}

\begin{proposition}
The action of $\brn$ on $\Xi$ by translation is isometric
and strongly continuous.
\end{proposition}

\begin{proof}
That the action is isometric is evident from the definition
of the norm. Since the action is isometric, it suffices to
show that the action is strongly continuous on $C_c(\brn)$,
since $C_c(\brn)$ is
dense by Proposition \ref{prop1c}. Now if $\xi \in C_c(\brn)$
then for $z$'s and $x$'s restricted
to given compact sets there is only a finite number of
$p$'s in $\bzn$ such that $(T_z\xi)(x-p) - \xi(x-p) \neq 0$.
The strong continuity then follow from the fact that $\xi$
is uniformly continuous since it is in $C_c(\brn)$.
\end{proof}

\begin{corollary}
The infinitely
differentiable functions of compact support are dense
in $\Xi$.
\end{corollary}

\begin{proof}
Since translation is strongly continuous, the usual
smoothing argument works, consisting of convolving functions
in $C_c(\brn)$ by smooth functions of compact support.
\end{proof}

By definition \cite{La}, a ``right
Hilbert $C^*$-module'' over a $C^*$-algebra $\ca$ is a right
$\ca$-module which is equipped with an $\ca$-valued inner
product and which is complete for the corresponding norm
(defined much as in Definition \ref{def1a}).
In terms of this definition we can summarize the early
part of the above discussion by:

\begin{proposition}
\label{prop1b}
When $\Xi$ is equipped with
its $\ca$-valued inner product and right action of $\ca$,
it is a right Hilbert
$C^*$-module over $\ca$ (or a ``right Hilbert $\ca$-module'').
\end{proposition}

We remark that many of the considerations above apply in the
setting of an arbitrary locally compact group $G$ (or a Lie
group) and a discrete
subgroup $G_0$ such that $G/G_0$ is compact.

For each $j$ we find that
$\ca_j = C(\brn/B^{-j}\bzn)$. Since
$B^{-j}\bzn \supset B^{-(j+1)}\bzn$, we have
$\ca_j \subset \ca_{j+1}$ as expected,
where again we view elements of
$C(\brn/B^{-j}\bzn)$ as functions on $\brn$. We also have
\[
(\a(f))(x) = f(Bx)
\]
for $f \in \ca_j$. For each algebra $\ca_j$ we have a corresponding
$\ca_j$-valued inner product on a dense subspace of
$L^2(\brn)$. But in fact these dense subspaces all
coincide with $\Xi$, because the corresponding norms are
all equivalent. In fact we have:

\begin{proposition}
For any given $j \in \bz$ let $\G$ be the image of $B^{-j}\bzn$ in
$\brn/ (B^{-(j+1)}\bzn)$, so that $\G$ is a finite group.
Then for any $\xi, \eta \in \Xi$ we have
\[
\<\xi, \eta\>_{\ca_j}(x) =
\sum_{s \in \G}\<\xi, \eta\>_{\ca_{j+1}}(x-s)
\]
for all $x \in \brn$.
We also have
\[
\|\xi\|_{\ca_{j+1}} \leq \|\xi\|_{\ca_j} \leq \d^{-1}\|\xi\|_{\ca_{j+1}}.
\]
\end{proposition}

\begin{proof}
   For simplicity
of notation we just compare $\ca = \ca_0$ with $\ca_1$,
so that $j=0$.
Pull back $\G$ to a set $C$ of coset
representatives for $B^{-1}\bzn$ in $\bzn$, so that
\[
|C| = |\det(B^{-1})| = |\det(A)| = \d^{-2},
\]
where $|C|$ denotes the number of elements in $C$.
Then for any $x \in \brn$ we have
\[
\<\xi, \eta\>_{\ca}(x) =
\sum_{\bzn} (\bar\xi\eta)(x - p)
=\sum_{c \in C} \sum_{q \in \bzn}
(\bar\xi\eta)(x - c - B^{-1}q) =
\sum_{s \in \G}\<\xi, \eta\>_{\ca_{1}}(x-s).
\]
The second inequality follows easily from this. The first
inequality follows immediately from $B^{-1}\bzn \subset \bzn$.
\end{proof}

\begin{corollary}
\label{cor1a}
Let $\xi, \eta \in \Xi$. If $\<\xi, \eta\>_{\ca_j} = 0$ for
some $j$, then $\<\xi, \eta\>_{\ca_k} = 0$ for all $k<j$.
\end{corollary}

We point out the following
view of the above discussion.
It is natural to define a
``conditional expectation'', $E$, from $\ca_1$ onto $\ca$
by averaging over $\G$. That is,
\[
E(g)(x) = \d^{2} \sum_{s \in \G} g(x-s)
\]
for $g \in \ca_1$. It is clear that
$\|E(g)\|_\infty \leq \|g\|_\infty$.
Then, the first part of the above proposition says that
\[
\<\xi, \eta\>_{\ca_j} = \d^{-2} E(\<\xi, \eta\>_{\ca_{j+1}}),
\]
while the second inequality of the above proposition amounts to
\[
\|\<\xi, \xi\>_{\ca_j}\|_\infty =
\d^{-2}\|E(\<\xi, \xi\>_{\ca_{j+1}})\|_\infty
\leq \d^{-2}\|\<\xi, \xi\>_{\ca_{j+1}}\|_\infty.
\]

The following formulas, whose versions in the time domain
were mentioned earlier, are verified by straightforward calculations.

\begin{proposition}
\label{prop1t}
For any $\xi, \eta \in \Xi$ and $f \in \ca_j$ we have
\[
D(\xi f) = (D\xi)(\a(f)) \quad\quad \mathrm{and}
\quad\quad \<D\xi, D\eta\>_{\ca_{j+1}} = \a(\<\xi, \eta\>_{\ca_j}) .
\]
\end{proposition}

\section{Projective modules and frames}

Suppose now that $\var_1,\dots , \var_m \in \Xi$ and
that $\<\var_i, \var_j\>_\ca = \d_{ij}1_\ca$.
Then it is easily checked that the transformation
\[
(a_1, \dots , a_m) \mapsto \sum^m \var_j a_j
\]
from $\ca^m$ into $\Xi$ is a right $\ca$-module
homomorphism which is isometric for the $\ca$-valued
inner products. Here $\ca^m$ is viewed as a right
$\ca$-module in the evident way, and is equipped with
the natural inner product
\[
\<(a_i), (b_j)\>_\ca = \sum^m a_k^* b_k.
\]
Then the $\ca$-submodule of $\Xi$ generated by
$\var_1, \dots ,\var_m$ is norm-closed, and is
isometrically isomorphic to $\ca^m$. Now the $\ca$-modules
which are isomorphic to ones
of the form $\ca^m$ are exactly what are called (finitely
generated) free $\ca$-modules. We observe that in the
time domain the $\ca$-module generated by the $\var_i$'s
is exactly the closed span of the translates of the $\var_i$'s by
elements of $\bzn$ (and these translates are all
orthogonal for the ordinary inner product on $L^2(\brn)$).
This is the standard situation considered
for wavelets and multi-wavelets.

The main thrust of our paper is that one can use somewhat
more general $\ca$-modules for the same purposes,
namely (finitely generated) projective modules. Their
usefulness in connection with wavelets was already
demonstrated in \cite{pacrief}. For the reader's convenience
we repeat here our brief discussion in \cite{pacrief} of the
definition of projective modules, and of a few facts
about them.

By definition,
a projective module is (isomorphic to) a direct summand of a
free module. In this paper all projective modules will be
taken to be finitely generated.
In the next section we will give specific examples of non-free
projective $\ca$-modules for the case $n=2$.
Let us view our free module as ${\mathcal A}^m$ for
some integer $m$, with ``standard basis'' $\{e_j\}$.
Then the definition
means that there is an $m \times m$ matrix $P$ with
entries in ${\mathcal A}$ which is a projection, that is
$P^2 = P$, such that our projective module $V$ is of
the form $V = P{\mathcal A}^m$.
In this (C*-algebraic) setting
a standard argument (see 5Bb in \cite{We}) shows that $P$ can
be adjusted so that it is also ``self-adjoint'' in the evident sense.
Let us then set $\xi_j = Pe_j$ for each $j$.
For any $\eta \in V$ we have
$$
\eta = P\eta = P(\sum e_j\langle e_j, \eta \rangle_{\mathcal A})
= \sum Pe_j \langle e_j, \eta \rangle_{\mathcal A}
\\ = \sum \xi_j \langle e_j, P\eta \rangle_{\mathcal A}
$$
$$
= \sum \xi_j \langle Pe_j, \eta \rangle_{\mathcal A}
= \sum \xi_j \langle \xi_j, \eta \rangle_{\mathcal A}   .
$$
There is no reason to expect that the $\xi_j$'s will be
independent over ${\mathcal A}$, much less orthonormal. But anyone
familiar with wavelets will feel comfortable about referring to
the $\xi_j$'s as a ``module frame'' for $V$.
In the more general setting of projective modules over C*-algebras
the above reconstruction formula
\begin{equation*}
\label{eq pmf}
         \eta\;=\;\sum \xi_j\;\langle \xi_j, \eta
\rangle_{\mathcal A}
\end{equation*}
appears, in different notation, already in \cite{Rie}.

In the study of Hilbert modules the concept of module frames
is extremely useful. We modify for our purposes some of the
definitions and results about module frames  given in the
paper of Larson and Frank \cite{LaFr98}, starting with the
following definition:
\begin{definition}
\label{def2framefinite}
Let $V$ be a finitely-generated right Hilbert $C^*$-module over a
unital $C^*$-algebra $\ca$. By a {\em standard module frame} for $V$
we mean a finite collection $\var_1,\dots , \var_m$ of
elements of $V$ for which the reconstruction formula
\[
v = \sum_j \var_j\<\var_j, v\>_\ca
\]
holds for all $v \in V$.
\end{definition}

For use in Theorem \ref{thm2x} we also give the definition of
standard module frames for Hilbert $C^{\ast}$-modules
that are not finitely generated, due to Frank and
Larson \cite{LaFr98, LaFr00}.

\begin{definition}
\label{def2frameinf}
Let $\Xi$ be a right Hilbert $C^*$-module over a
unital $C^*$-algebra $\ca$. We say that a countable subset
$\{\var_j\}_{j\in I}$ of $\Xi$ is a {\em standard module
frame} for $\Xi$
if for every $\xi\in \Xi,$
\[
\langle \xi, \xi \rangle_{\ca}\; = \; \sum_{j\in I}
\langle \xi,  \var_j \rangle_{\ca} \langle \var_j, \xi \rangle_{\ca},
\]
where the sum on the right-hand side converges in norm in $\ca.$
\end{definition}

Frank and Larson have shown in Theorem 4.1 of 
\cite{LaFr00} that if  $\{\var_j\}_{j\in I}$
is a standard module frame for $\Xi$,
then the reconstruction formula
\begin{equation}
\label{eq2recform}
\xi \;=\;\sum_{j\in I} \var_j\<\var_j, \xi\>_\ca
\end{equation}
holds for every $\xi\in \Xi,$ where the sum converges in norm.
In the case where the index set $I$ is finite,
one can easily verify that
the two definitions \ref{def2framefinite} and \ref{def2frameinf}
of standard module frames are equivalent. For the moment we concentrate on Hilbert $C^{\ast}$  submodules of $\Xi$ having standard module frames of finite cardinality; examples of standard module frames of countably infinite cardinality will be given in Theorem \ref{thm2x} in the next section. 

The discussion before the definitions
shows that any (finitely generated)
projective module has a
standard module frame, in the sense of Definition \ref{def2framefinite}.
We now show, in two steps, that the converse is true.

\begin{proposition}
\label{prop2a}
Let $V$ be a submodule of some Hilbert $C^*$-module $\Xi$ over
some unital $C^*$-algebra $\ca$ which has a standard module frame $\{\var_j\}$ of finite cardinality. Then $V$
is an $\ca$-module direct summand of $\Xi$; a projection
operator, $P$, from $\Xi$ onto $V$ is given by
\[
P(\xi) = \sum \var_j\<\var_j, \xi\>_\ca .
\]
\end{proposition}

\begin{proof}
By Definition \ref{def2framefinite} of a standard module
frame it is clear that the
restriction of $P$ to $V$ is the identity operator on $V$.
Thus $P^2 = P$. A simple calculation shows that $P^* = P$ in
the sense that $\<P\xi, \eta\>_\ca = \<\xi, P\eta\>_\ca$.
\end{proof}

We remark that the above proposition is false, in general, for submodules $V$ of $\Xi$ having a countably infinite standard module frame $\{\var_j\}_{j=1}^{\infty}$.  The problem arises because the sum $\sum \var_j\<\var_j, \xi\>_\ca$ need not converge in norm for every $\xi$ in $\Xi,$ although the sum will converge in norm for $\xi\in V.$  We thank the referee for bringing this point to our attention. 

\begin{proposition}
\label{prop2fr}
Let $V$ be a Hilbert $C^*$-module
over $\ca$ which has a standard module frame of finite cardinality.
Then $V$ can be embedded into some free module $\ca^m$, with
preservation of the $\ca$-valued inner products. Furthermore,
$V$ is a projective $\ca$ module.
\end{proposition}

\begin{proof} Let $\var_1, \dots , \var_m$ be a standard module
frame for $V$. Define $T: V \mapsto \ca^m$ by
\[
Tv = (\<\var_1, v\>_\ca, \dots, \<\var_m, v\>_\ca).
\]
Simple calculations show that $T$ gives the desired embedding.
The fact that $V$ is projective then follows from
Proposition \ref{prop2a}.
\end{proof}

We now describe the relationship with ordinary frames for $L^2$-spaces
within our specific setting of $\ca = C(\btn)$ and
our $\Xi \subset L^2(\brn)$.
For $q \in \bzn$ we let $e_q \in C(\btn)$ be defined by
$e_q(x) = \exp(2\pi i x\cdot q)$.

\begin{proposition}
\label{prop2t}
Let $V$ be a projective $\ca$-submodule of $\Xi$, and let
$\{\var_k\}$ be a standard module frame for $V$. Then $\{\var_k e_q\}$
is a normalized tight frame for the closure of $V$ in $L^2(\brn)$.
\end{proposition}

\begin{proof}
For any $\xi, \eta \in \Xi$ we have, by easily justified
manipulations,
\begin{eqnarray*}
\lefteqn{\sum_{k, q}\overline{\<\var_k e_q, \xi\>}_{L^2}
\<\var_k e_q, \eta\>_{L^2}
= \sum_{k, q} \int \overline{\<\var_k e_q, \xi\>}_\ca
\int \<\var_k e_q, \eta\> _\ca} \\
&=&\sum_k\big(\sum_q\int e_q\overline{\<\var_k, \xi\>}_\ca
\int \bar e_q\<\var_k, \eta\>_\ca\big)
= \sum_k\int \overline{\<\var_k, \xi\>}_\ca\<\var_k, \eta\>_\ca \\
&=& \int\<\sum_k \var_k\<\var_k, \xi\>_\ca, \eta\>_\ca
= \int\<\xi, \eta\>_\ca = \<\xi, \eta\>_{L^2}.
\end{eqnarray*}
This relation then extends to the closure of $V$ in $L^2(\brn)$.
Thus the basic condition in the definition of a normalized
tight frame has been verified. A similar calculation shows that
\[
\xi = \sum_{k,q} (\var_k e_q)\<(\var_k e_q), \xi\>_{L^2}
\]
for any $\xi \in V$, and so for any $\xi$ in the closure of $V$ in $L^2(\brn)$.
\end{proof}

In the same way, standard module frames for projective
submodules of $\Xi$ over any
$\ca_j$ give normalized tight frames for their
closures in $L^2(\brn)$.

\section{Projective multiresolution analyses}

The main point of our paper is that it is possible to construct
multiresolution analyses from projective modules which are not
free. Explicit constructions of such multiresolution
analyses will be given in the next sections. Here we just
give the definition of such analyses, in terms
of the notation which we have been using, and then explore
some of their general properties. \begin{definition}
\label{defpr}
By a {\em projective multiresolution analysis} for dilation by
$A$ we mean a family $\{V_j\}_{j \in \bz}$ of subspaces of
$\Xi$ such that
\begin{enumerate}
\item $V_0$ is a projective $\ca$-submodule of $\Xi$.
\item $V_j = D^j(V_0)$ for all $j$.
\item $V_j \supset V_{j-1}$ for all $j$.
\item $\bigcup^\infty V_j$ is dense in $\Xi$
\item $\bigcap^\infty_{-\infty} V_j = \{0\}$.
\end{enumerate}
\end{definition}

We comment on the relationship between our definition and the definition of projective multiresolution analyses given by Zimmermann in his thesis \cite{Z1}. In Chapter VII of \cite{Z1}, a projective multiresolution analysis with scale $k\in\mathbb N\backslash\{1\}$ is defined to be a bisequence of continuous projections $\{P_m\}_{m\in\mathbb Z}$ defined on a translation-invariant locally convex topological vector space $X(\mathbb R)$ of functions or distributions on the real line which satisfy:
\begin{enumerate}
\item [[P1]] $P_{m_1}P_{m_2} = P_{\text{min}\{m_1,m_2\}};$
\item [[P2]] For every $f\in X(\mathbb R),\;\lim_{m\rightarrow -\infty}P_m(f)= 0$ and $\lim_{m\rightarrow\infty}P_m(f)=f;$
\item [[P3]] For every $m\in\mathbb Z,\;D_kP_m = P_{m+1}D_k,$ where $D_k$ is the operator on $X(\mathbb R)$ of dilation by $k;$ 
\item [[P4]] For all $m,n\in \mathbb Z,\;\tau_{k^{-m}n} P_m = P_m\tau_{k^{-m}n},$ where for $t\in\mathbb R,\;\tau_t$ is the operator of translation by $t$ defined on $X(\mathbb R),$ and $k\in\mathbb N\backslash\{1\}$ is the dilation factor;
\item [[P5]] There exists $\phi\in X(\mathbb R),\;\phi^{\ast}\in X(\mathbb R)^{\ast}$ such that 
$\{\tau_n\phi,\tau_n\phi^{\ast}\}_{n\in\mathbb Z}$ is a biorthogonal system with 
$P_0(f)=\sum_{n\in\mathbb Z} \langle f,\tau_n\phi^{\ast}\rangle \tau_n\phi.$
\end{enumerate}
We note that given a sequence of subspaces $\{V_j\}_{j \in \bz}$ of $\Xi$ satisfying Definition \ref{defpr} it is possible to form continuous projections $\{P_m\}$ on $\Xi$ that satisfy conditions analogous to [P1]-[P4] of Zimmermann's definition with respect to the dilation $A$, first for the case $m\geq 0$ by using the constructions of projections given in Proposition \ref{prop2a}, and then modifying this method for the case $m<0$.  However we have no analogue of [P5] in our definition, and our main aim is to construct $C^{\ast}$-module frames, as opposed to Zimmermann's emphasis on constructing an analogue of biorthogonal frames in his setting.  
We thank the referee for bringing Zimmermann's work to our attention. 

Returning to Definition \ref{defpr}, we quickly point out that condition (5) above is
actually redundant in our situation. That is, even more
generally, assuming only (1) and (2) we obtain the following result:

\begin{proposition}
\label{prop2i}
Let $V$ be any projective $\ca$-submodule of \ $\Xi$,
and set $V_j = D(V_{j-1})$ for all $j$. Then
$\bigcap^\infty_{-\infty} V_j = \{0\}$.
\end{proposition}

\begin{proof}
Let $\xi \in \bigcap^\infty_{-\infty} V_j$. Then $D^{-j}\xi \in V$
for each $j$. Let $\{\eta_k\}$ be a frame for $V$. Thus
$D^{-j}\xi = \sum_k\eta_k\<\eta_k, D^{-j}\xi\>_\ca$, so that
\[
\xi = \sum_k D^j(\eta_k\<\eta_k, D^{-j}\xi\>_\ca)
\]
for each $j$. We will try to show that $\<g, \xi\> = 0$ for every
$g \in C_c(\brn)$. Now for such $g$ we have
\[
\<g, \xi\> =
\sum_k \<g, D^j(\eta_k\<\eta_k, D^{-j}\xi\>_\ca)\>.
\]
So we examine the individual terms of this sum.

For fixed $k$ let $F_j = \<\eta_k, D^{-j}\xi\>_\ca$. Then, using
the Cauchy-Schwarz inequality of Proposition \ref{prop1cs},
and then Proposition \ref{prop1a}, we have
\begin{eqnarray*}
\|F_j\|^2_2 &=& \int_{\btn} |F_j|^2 =
\int |\<\eta_k, D^{-j}\xi\>_\ca|^2  \\
&\leq& \int \<\eta_k, \eta_k\>_\ca\<D^{-j}\xi, D^{-j}\xi\>_\ca
\leq \|\<\eta_k, \eta_k\>_\ca\|_\infty
\int \<D^{-j}\xi, D^{-j}\xi\>_\ca \\
&=& \|\eta_k\|_\ca^2 \|D^{-j}\xi\|^2_2 = \|\eta_k\|^2_\ca \|\xi\|^2_2.
\end{eqnarray*}
Thus
\[
\|F_j\|_2 \leq  \|\eta_k\|_\ca \|\xi\|_2.
\]
Note that the right-hand side is independent of $j$.

Because $A$ is a dilation, a bit of examination shows that we
can hope that as $j \rightarrow -\infty$ the ``mass'' of
$ D^j(\eta_k\<\eta_k, D^{-j}\xi\>_\ca)$ will converge towards
$0 \in \brn$. In order to quantify this idea, for any $r>0$
we define the function $\chi_r$ on $\brn$ by giving it value
$1$ outside the ball of radius r, and value $0$ inside the ball.

\begin{lemma}
Let $\eta \in \Xi$. For any $\e>0$ there is an $r>0$ such that
for all $F \in \ca$ we have
\[
\|\chi_r\eta F\|_2 \leq \e\|F\|_2.
\]

\begin{proof}
The sum $\sum |\eta(x-p)|^2$ converges uniformly for $x \in I^n$,
and so we can find a finite subset, $S$, of $\bzn$ such that
$\sum_{p \notin S} |\eta(x-p)|^2 < \e^2$ for all $x \in I^n$.
If $p \in S$ and $x \in I^n$ then $\|x-p\| \leq \|x\| + \|p\|
\leq \sqrt{n} + \|p\|$. Set
\[
r = \sqrt{n} + \max\{\|p\|: p \in S\} + 1.
\]
Thus if $\chi_r(x-p) \neq 0$ for an $x \in I^n$,
then $\|x-p\| \geq r$, so that $p \notin S$.
Consequently
\begin{eqnarray*}
\|\chi_r\eta F\|^2_2 &=& \int_{\brn} \chi_r |\eta F|^2
= \sum_p \int _{p+I^n} \chi_r |\eta F|^2
= \sum_p \int_{I^n} \chi_r(x-p) |(\eta F)(x-p)|^2 \ dx  \\
&\leq& \sum_{p \notin S} \int_{I^n} |\eta(x-p)|^2 |F(x)|^2 \ dx
= \int_{I^n} \sum_{p \notin S}|\eta(x-p)|^2 |F(x)|^2 \ dx
\leq \e^2 \|F\|^2_2.
\end{eqnarray*}
\end{proof}
\end{lemma}

We keep the notation of the Lemma.
Since $D$ is unitary, it follows from the Lemma that
$\|D^j(\chi_r \eta F)\|_2 \leq \e\|F\|_2$.
Set $\tilde \chi_r = 1 - \chi_r$. Then
\begin{eqnarray*}
|\<g, D^j(\eta F)\>| &\leq& |\<g, D^j(\chi_r\eta F)\>|
+ |\<g, D^j(\tilde \chi_r\eta F)\>| \\
&\leq& \e \|g\|_2\|F\|_2 + |\<g, D^j(\tilde \chi_r\eta F)\>|.
\end{eqnarray*}
Suppose now that $0$ is not in the (closed) support of $g$.
Notice that $(D^j(\tilde\chi_r\eta F))(x) =
\d^j ((\tilde \chi_r \eta F)(B^jx))$. If
$\tilde \chi_r(B^jx) \neq 0$ then $\|B^jx\| \leq r$, so
that $\|x\| \leq \|B^{-j}\|r$. Since $A$ is a dilation,
the spectral radius of $B$ is $<1$. It follows from the
spectral radius formula that $\|B^{-j}\| \rightarrow 0$
as $j \rightarrow -\infty$. Thus for sufficiently negative $j$
the supports of $g$ and $\tilde \chi_r \eta F$ are disjoint, so that
\[
|\<g, D^j(\eta F)\>| \leq \e\|g\|_2\|F\|_2.
\]
Putting all this together with our calculations before the Lemma,
we find that for $j$ sufficiently negative so that the above
estimates apply for all of our finite number of indices $k$, we have
\[
|\<g,\xi\>| \leq \e(\sum_k\|g\|_2 \|\eta_k\|_\ca \|\xi\|_2).
\]
Since $\e$ is arbitrary, it follows that $\<g, \xi\> = 0$. Thus
we have shown that $\xi$ is orthogonal to all $g \in C_c(\brn)$
whose support does not contain $0$. This is sufficient to
conclude that $\xi = 0$, as desired.
\end{proof}

We now show that condition (4) of Definition \ref{defpr} holds
under a natural condition which is closely related to the
familiar condition on scaling functions that $\hat\var(0) = 1$.

\begin{proposition}
\label{prop2d}
Let $V$ be a projective $\ca$-submodule of $\Xi$ which satisfies
conditions (1), (2) and (3) of Definition \ref{defpr}. If there
is at least one
$\xi \in V$ such that $\xi(0) \neq 0$, then $V$ satisfies
condition (4) of Definition \ref{defpr}.
\end{proposition}

\begin{proof}
The direction of our proof is suggested by the proof of theorem
1.7 of chapter 2 of \cite{HW}. Assume that $V$ satisfies the
conditions of the proposition, and set
$\cv = \bigcup^\infty_{-\infty} V_j$.  Set
$\cb = \bigcup^\infty_{-\infty} \ca_j$, so that $\cb$ is a
unital $*$-subalgebra of $C_b(\brn)$ which separates the
points of $\brn$. Let $\overline{\cv}$ denote the closure of $\cv$
in $\Xi$. Notice that $\cv$ is a $\cb$-module
since $V_j$ is an $\ca_j$-module for
each $j$.

\begin{lemma}
\label{lem2p}
The subspace $\overline{\cv}$ is closed under pointwise multiplication
by $C_b(\brn)$.
\end{lemma}

\begin{proof}
Let $\xi \in \overline{\cv}$ and $F \in C_b(\brn)$.
By Proposition \ref{prop1c} we know that $\xi F \in \Xi$, and so
it suffices to show that $\xi F$ can be approximated in the
norm of $\Xi$ by elements of $\cv$. Let $\e >0$ be given. Let
$S$ be a finite subset of $\bzn$ such that for any $x \in I^n$
we have
$\sum_{p \notin S}|\xi(x-p)|^2 < (\e/\|F\|_\infty)^2$.
   Let $\overline{\cb}$ denote the closure
of $\cb$ in $C_b(\brn)$. Then $\overline{\cb} \cong C(Y)$ for some
compact space $Y$, and there is a natural continuous
injection of $\brn$ into $Y$ as a dense subset (not open --- we
are close to the Bohr compactification). Since
$I^n - S$ is a compact subset of $\brn$, its image
is a compact subset of $Y$, on which the restriction to
$I^n -S$ of $F$ can be considered to be
a continuous function. Thus by the Tietze extension
theorem this function extends to a continuous function
on $Y$, and so we can choose a $G \in \cb$ such that
$\|G\|_\infty \leq \|F\|_\infty$ and
$|F(x) - G(x)| < \e/\|\xi\|_\ca$ for every $x \in I^n-S$.
Set $H=F-G$. Then for any $x \in I^n$ we have
\begin{eqnarray*}
\<\xi H, \xi H\>_\ca(x) &=&
\sum_{p \in S}|(\xi H)(x-p)|^2 + \sum_{p \notin S}|(\xi H)(x-p)|^2 \\
&\leq& (\e/\|\xi\|_\ca)^2 \sum_{p \in S}|\xi(x-p)|^2
   +2\|F\|_\infty^2 \sum_{p \notin S}|\xi(x-p)|^2  \\
&\leq& \e^2+ 2\e^2 = 3\e^2.
\end{eqnarray*}
It follows that $\|\xi F - \xi G\|_\ca < \sqrt{3} \e$.
\end{proof}

We now return to the proof of Proposition \ref{prop2d}. We show
that $C_c(\brn) \subset \overline{\cv}$. The conclusion of the
proposition then follows from Proposition \ref{prop1c}. So
let $\eta \in C_c(\brn)$, and let $K$ denote the support of $\eta$.
Choose a $\xi\in V$ such that $\xi(0) \neq 0$. Since $\xi$ is
continuous, $\xi(x) \neq 0$ for all $x$ in some neighborhood
of $0$. As seen in the proof of Proposition \ref{prop2i},
we can choose a
sufficiently large $j$ that $\|B^j\|$ is very small, small
enough that $\xi \circ (B^j)$ is bounded away
from $0$ on $K$ . Note
that $\xi \circ (B^j) \in V_j$ by definition.
Set $g = \eta(\xi \circ B^j)^{-1}$.
Then $g \in C_c(\brn)$. Since $\eta = (\xi\circ B^j)g$,
it follows from Lemma \ref{lem2p}
that $\eta \in \overline{\cv}$ as desired.
\end{proof}

If $k < j$ then $\ca_k \subset \ca_j$, and so
we can also view $V_j$ as an $\ca_k$-module.

\begin{proposition}
\label{prop2k}
If $k < j$, than
$V_j$ is projective as an $\ca_k$-module.
\end{proposition}
\begin{proof} Let $V$ be any projective $\ca_j$-module. Then
there is an $\ca_j$-module $W$ such that
$V\oplus W \cong (\ca_j)^m$ for some integer $m$.
If $V_{\ca_k}$ denotes $V$ viewed as an $\ca_k$-module, and
similarly for other $\ca_j$-modules, then we have
\[
V_{\ca_k} \oplus W_{\ca_k} \cong ((\ca_j)^m)_{\ca_k}
\cong ((\ca_j)_{\ca_k})^m.
\]
We can thus obtain the desired conclusion if we know that
$A_j$ is projective as a right $\ca_k$-module. But $\ca_j$
is actually free as a right $\ca_k$-module, with the proof
being essentially proposition 1.1 of \cite{pacrief}. Indeed,
if we view matters in the time domain, then any set of coset
representatives for $A^{-k}\bzn $ in $A^{-j}\bzn$ will
provide a module basis.
\end{proof}

We remark that in the setting of the above proposition a standard
module frame for $V_j$ as $\ca_k$-module can be constructed as
follows from a standard module frame $\var_1, \dots , \var_m$
for $V_j$ as an $\ca_j$-module. For ease of notation we assume that
$k = 0$, so that $j>0$. We note first that
there is a natural $\ca$-valued inner product on $\ca_j$.
The image of $\bzn$ in $\brn/B^{-j}\bzn$ is a finite
group, say $\G$, which we can pull back to a set, $C$, of coset
representatives for $B^{-j}\bzn$ in $\bzn$. Define a linear map $E$
from $\ca_j$ onto $\ca$ by $E(f)(x) = \sum_{p \in C}f(x-p)$. (Then
$E$ is a multiple of a ``conditional expectation'' as discussed
after Corollary \ref{cor1a}.)
We define
our $\ca$-valued inner product on $\ca_j$ by
$\<f, g\>_\ca = E(\bar fg)$. The context must distinguish this inner
product from that on $\Xi$.

Since $\ca_j$ is free as an $\ca$-module, we can find an  orthonormal
module basis $b_1, \dots , b_q$ for $\ca_j$ as an
$\ca$-module. Then
$\{\var_i b_k: 1 \leq i \leq m, 1 \leq k \leq q\}$ will be a standard
module frame for $V_j$ as an $\ca$-module. To see this,
let $\xi \in V_j$. Then
\begin{eqnarray*}
\xi &=& \sum \var_i\<\var_i, \xi\>_{\ca_j} =
\sum \var_i (\sum b_k\<b_k, \<\var_i, \xi\>_{\ca_j}\>_\ca)
= \sum(\var_i b_k)\<b_k, \<\var_i, \xi\>_{\ca_j}\>_\ca  \\
&=& \sum(\var_i b_k)E(\bar b_k\<\var_i, \xi\>_{\ca_j})
= \sum (\var_i b_k)E(\<\var_i b_k, \xi\>_{\ca_j})
= \sum (\var_i b_k)\<(\var_i b_k), \xi\>_\ca ,
\end{eqnarray*}
so the reconstruction formula holds, as desired.

Suppose now that $\{V_j\}$ is a projective multiresolution analysis.
Then $V_1$ is a projective $\ca_1$-module, which contains
$V = V_0$, but $V$ is not an $\ca_1$-submodule. However, from
Proposition \ref{prop2k} we can consider $V_1$ to be a projective
$\ca$-module, and then $V$ is a projective $\ca$-submodule of $V_1$.
{}From Proposition \ref{prop2a} it follows that $V$ has a unique
``orthogonal complement'', using the $\ca$-valued inner product.
Let us denote this complement by $W_0$. Then $W_0$ is again
a projective $\ca$-module, and so it will have a standard module frame.
These are our wavelets. One way of obtaining
a standard module frame for $W_0$ is as follows. Choose
a standard $\ca$-module frame for $V$, and apply $D$ to it to obtain
a standard $\ca_1$-module frame for $V_1$. Then apply
Proposition \ref{prop2k} and the discussion which
follows it to obtain a standard $\ca$-module frame, say
$\var_1, \dots , \var_m$, for $V_1$. For the chosen standard
$\ca$-module frame for $V$ use Proposition \ref{prop2a} to construct
the orthogonal projection $P$ from $V_1$ onto $V$. Then $Q = P - I$ is the
orthogonal projection of $V_1$ onto $W_0$. The discussion before
Definition \ref{def2framefinite} shows that $Q\var_1, \dots , Q\var_m$
will then form a standard $\ca$-module frame for $W_0$, hence they will
be our wavelets.

  From Proposition \ref{prop2t} we see that the $\psi_k e_q$'s will form
an ordinary normalized tight frame for the closure of $W_0$
in $L^2(\brn)$. In the same way we can let $W_j$ be the
orthogonal complement of $V_j$ in $V_{j+1}$ as $\ca_j$-modules.
As is usual for multiresolution analyses, we have:

\begin{proposition}
\label{prop2w}
The operator $D$
will carry $W_j$ onto $W_{j+1}$, and will carry a standard module frame
for $W_j$ as $\ca_j$-module to one for $W_{j+1}$ as $\ca_{j+1}$-module,
with a similar relation for
the corresponding normalized tight frames.
\end{proposition}

\begin{proof}
This follows from several straight-forward calculations
using Proposition \ref{prop1t} together with the fact that $D$
is a unitary operator on $L^2(\brn)$.
\end{proof}

From
Corollary \ref{cor1a}, condition (4) of Definition \ref{defpr},
and Proposition \ref{prop2w}, we obtain:

\begin{theorem}
\label{thm2w}
The $\ca$-submodules $V$ and $W_j$ for $j \geq 0$ are all mutually
orthogonal for the
$\ca$-valued inner product, and the algebraic
sum \ $V\oplus_{j\geq 0}^{\infty} W_j$ \
is dense in $\Xi$.
\end{theorem}

Let $\overline V$ and ${\overline W}_j$ denote the closures of $V$ and $W_j$
in $L^2(\brn)$. Since $\Xi$ is dense in $L^2(\brn)$ it follows
from Proposition
\ref{prop1a} that:

\begin{corollary}
With notation as above we have the Hilbert-space direct sums
\[
L^2(\brn) =  {\overline V}\oplus_{j\geq 0}^{\infty} {\overline W}_j
= \oplus_{j \in \bz} {\overline W}_j  .
\]
\end{corollary}

Of course in view of this, the union of a normalized tight frame for
$\overline V$
together with ones for each
of the ${\overline W}_j$'s for $j > 0$ is a normalized
tight frame for $L^2(\brn)$. In the same way any collection of
normalized tight frames for each of the ${\overline W}_j$'s
as $j$ runs over $\bz$ will be a
normalized tight frame for $L^2(\brn)$. The usual pattern
for multi-wavelets would be to obtain the tight frames for each
${\overline W}_j$ by first finding a Hilbert $\ca_0$-module
basis for $W_0$ (as well as for $V$), then applying Proposition
\ref{prop2w}, and then applying the $j$-versions
of Proposition \ref{prop2t}.

  From the discussion of topological direct sums of Hilbert
$C^*$-modules given in Chapter 1 of \cite{La} we find
that the topological sum $V\oplus_{j\geq 0}^{\infty} W_j$
within $\Xi$ can be identified with the set
$$\{v+\{w_j\}_{j\geq 0}\, : w_j \in W_j \; \;
\text{and}\; \;\sum_{j\geq 0}\langle w_j, w_j
\rangle_{\ca}\; \; \text{is convergent in}\;\ca\},$$
which is complete.  Thus the topological
sum $V\oplus_{j\geq 0}^{\infty} W_j$
is a complete submodule of $\Xi$ which
contains $V_j$ for all$j\geq 0.$ It follows by
condition $(4)$ of Definition \ref{defpr} that as a topological
sum we have $V\oplus_{j\geq 0}^{\infty} W_j = \Xi$ .

We have a similar result for standard module frames. Here we must view
each $W_j$ as an $\ca$-module, as permitted by Proposition
\ref{prop2k}. (This is not the usual point of view in
wavelet theory.) We use the definition of standard module frames
for non-finitely-generated modules
given in Definition \ref{def2frameinf} .

\begin{theorem}
\label{thm2x}
Let $V$ and $W_j$ for $j \geq 0$ be as in the
statement of Theorem \ref{thm2w}.
Then the union of a standard $\ca$-module frame for $V$ and
any collection of  standard $\ca$-module
frames for each of the $W_j$'s for
$j\geq 0,$ is a standard $\ca$-module frame for $\Xi.$
\end{theorem}

\begin{proof}
Let $\{\var_1, \dots , \var_m\}$ be a standard module
frame for $V$, and let
$\{\psi_{j,1},\cdots, \psi_{j,d_j}\}$ be a standard
$\ca$-module frame $W_j$ for each  $j\geq 0$. According to
Definition \ref{def2frameinf}
we must show that for every $\xi\in\Xi$ we have
$$
\langle \xi,\xi\rangle_{\ca}=\sum_{i=1}^m\langle
\var_i, \xi \rangle_{\ca} \langle \xi,  \var_i
\rangle_{\ca}+\sum_{j=0}^{\infty}\sum_{k=1}^{d_j}\langle
\xi, \psi_{j,k} \rangle_{\ca} \langle  \psi_{j,k}, \xi
\rangle_{\ca},
$$
with the sum on the right-hand side converging in norm in $\ca.$
Fixing $\xi\in\Xi,$ we use Theorem \ref{thm2w} to write
$$
\xi\;=\; v+\sum_{j=0}^{\infty}w_j,
$$
where $v\in V,\;w_j\in W_j$ for $j\geq 0,$ and
by Chapter 1 of \cite{La},
$$
\langle \xi, \xi \rangle_{\ca}\;=\;\langle v, v
\rangle_{\ca}+\sum_{j=0}^{\infty}\langle w_j, w_j \rangle_{\ca},
$$
where the sum on the right-hand side converges in norm in ${\ca}$.
But 
$$\langle v, v \rangle_{\ca}\;=\;\sum_{i=1}^m\langle
\var_i, v \rangle_{\ca} \langle v,  \var_i \rangle_{\ca},$$
while  
$$\langle w_j, w_j \rangle_{\ca}
\;=\;\sum_{k=1}^{d_j}\langle \xi_{j,k}, w_j\rangle_{\ca} \langle w_j,
\xi_{j,k} \rangle_{\ca}\;\;\;\text{for each}\; j\geq 0,$$ so that
$$
\langle \xi,\xi\rangle_{\ca}= \sum_{i=1}^m \langle v,
\var_i \rangle_{\ca} \langle \var_i, v
\rangle_{\ca}+\sum_{j=0}^{\infty}\sum_{k=1}^{d_j} \langle w_j,
\psi_{j,k} \rangle_{\ca} \langle \psi_{j,k}, w_j\rangle_{\ca},
$$
where again the right-hand side converges in norm in $\ca.$
Finally, we note that by the orthogonality of the modules involved,
$$
\langle \var_i, v \rangle_{\ca}\;=\;\langle \var_i, \xi \rangle_{\ca},
\qquad \text{and} \qquad
\langle \psi_{j,k}, w_j \rangle_{\ca}
\;=\;\langle \psi_{j,k}, \xi \rangle_{\ca},
$$
for all $1\leq i\leq m,\;j\geq 0,$ and $1\;\leq k\;\leq d_j.$
  From this we see that
the condition of Definition \ref{def2frameinf} is satisfied.
\end{proof}

\section{A review of finitely generated projective
$C(\mathbb T^2)$ modules}

In this section, we review results
on finitely generated projective modules over $\ca = C(\mathbb
T^2)$ which we will use often in later parts of the
paper. These results are module versions of known results
about complex vector bundles over $\bt^2$ which follow
easily from elementary vector bundle considerations, such
as those given in \cite{At} \cite{Hu}. The modules we consider
were constructed in  Section 3 of \cite{Rie1} to study the
rational rotation algebras ${\mathcal A}_{p/q}$ with
phase factor $e(p/q).$ The approach here will be
slightly different from that in \cite{Rie1}.
Throughout this section, unless otherwise specified, we again
view functions on $\mathbb T^n$ as functions defined on
$\mathbb R^n$ which are periodic modulo $\mathbb Z^n.$
The following proposition summarizes the results
from theorem 3.9 of \cite{Rie1}.

\begin{proposition}
For $q, \, a\in\mathbb Z$ with $q\neq 0$ let
$X(q,a)$ denote
the right $\ca$-module consisting of the space
of
continuous complex-valued functions $F$ on $\mathbb T
\times \mathbb
R$ which satisfy
$$
F(s,t-q)\;=\;e(as)F(s,t),
$$
with right $\ca$-module action given by
$$
(Ff)(s,t)\;=\;F(s,t)f(s,t),
$$
for $F\;\in\;X(q,a)$ and $f\;\in\;\ca.$
Then $X(q,a)$ is a finitely generated, projective
$\ca$-module.  The set
$\{X(q,a):\;q, \, a\;\in\;\mathbb Z,\; q>0\}$
parametrizes the
isomorphism classes of finitely generated projective
$\ca$-modules, in the sense that if $X$ is a finitely
generated
projective $\ca$-module, then there exist unique
values of $q$
and $a$ such that $X\;\cong\;X(q,a).$ For $q>0$ we
say that
$X(q,a)$ has dimension $q$ and twist $-a.$

We define an
$\ca$-valued inner product on $X(q,a)$ by
$$\langle F,\;G\rangle_{\ca}\;=\;\sum_{0\leq k\leq
q-1}{\overline {F(s,t-k)}}G(s,t-k),
$$
for $F,\, G\;\in\;X(q,a).$ The module $X(q,a)$ is complete for
the norm determined by this inner product, so that it is a
Hilbert $\ca$-module.
\end{proposition}

A simple calculation shows that $X(-q, a) \cong X(q,-a)$.
We note that there is an asymmetry in the treatment of
the variables in the formula defining the functions
in $X(q,a)$. This is just a matter of convention. But
it will affect later formulas which we use.

We now give a more detailed description of the $X(q,a)$'s.
First we fix some notation which will also be useful later.
For $q,\;a\,\in\,\mathbb Z$ with $q>0$,
and for $\beta \in \br$ define the function
$J_{q,a,\beta}:\;\mathbb
R^2\;\rightarrow\;\mathbb C$ by
$$
J_{q,a,\beta}(s,t)\;=\;e(nas) \quad \mathrm{for} \quad
t\;\in\;[\beta+nq,(\beta+nq)+q),\quad n\;\in\;\mathbb Z.
$$
We note that since $\mathbb R$ is the disjoint union
$\cup_{n\in\mathbb Z}[\beta+nq,(\beta+nq)+q),$ the above
formula
defines $J_{q,a,\beta}(s,t)$
on all of $\mathbb R^2.$  By construction,
$$
J_{q,a,\beta}(s,t-q)\;=\;e(as)J_{q,a,\beta}(s,t);
$$
but $J_{q,a,\beta}(s,t)$ is not an element of
$X(q,a)$ because
it is not continuous. However, if $f:\mathbb
R\;\rightarrow\;\mathbb
C$ is any $q$-periodic continuous function such that
$f(\beta)=0$
(so that $f(\beta+nq)=0$ for all $n\;\in\;\mathbb Z$ by
the $q$-periodicity of $f$), then it is easy to verify that
the function
$J_{q,a,\beta}(s,t)f(t)$ will be an element of $X(q,a).$
In fact, many infinitely
differentiable functions can be expressed in this way.

It is possible to show that if we fix
$\beta_0$ with $0<\beta_0 < q$, then the functions of the form
$J_{q,a,\beta}(s,t)f(t)$ for $\beta = 0$
together with those for $\beta = \beta_0$ , and for $f$
continuous, $q$-periodic and with $f(\beta)=0$,
form a set of $C(\mathbb T^2)$-generators for $X(q,a)$.
But as we do
not need such a result here, we do not include a
proof. The functions
$J_{q,a,\beta}$ together with low-pass filter functions
can be used to construct standard module frames for
each $X(q,a)$, and thus also the corresponding embeddings
into free modules. We give the procedure for $q=1$ and leave
the details for the general case
to the reader, since we will not need them later.
So let $q=1$ and let $a\,\in\,\mathbb Z$ be arbitrary.
Let $m_0$ be a continuous low-pass filter defined on
$\mathbb R$ for dilation by $2,$ so that, aside from
a constant such as $\sqrt{2}$, we have
$$
m_0(0)\;=\; 1,\;\text{and}
$$
$$
|m_0(t)|^2 + |m_0(t + 1/2)|^2 \;=\;1,\quad\text{so
that}\;m_0(1/2)=0.
$$
Set
$h_1(s,t)\;= m_0(t)J_{1,a,0}(s,t)$ and
$h_2(s,t)=m_0(t+ 1/2)J_{1,\;a,\; 1/2}(s,t).$
Then $h_1,\;h_2\;\in\;X(q,a),$
and it is easily verified that they form a standard module frame
for $X(1,\; a)$.
We refer the reader to
Proposition 2.1 of \cite{pac1} and to \cite{Rie1} for
ideas on how to go about extending this procedure for
general $q.$

\begin{theorem}
\label{thmrief}
(c.f. \cite{Rie1}, Theorem 3.9) Let
$q_1,\;q_2, \; a_1,\;a_2\;\in\;\mathbb Z$ with $q_1>0, \; q_2>0$
be given.
Then
$$X(q_1,a_1)\;\oplus\;X(q_2,a_2)\;\cong\;X(q_1+q_2,a_1+a_
2)$$ as
finitely generated projective $\ca$-modules.
In particular, cancellation holds for all finitely
generated $\ca$-modules; that is, if
$$X(q_1,a_1)\;\oplus\;X(q_2,a_2)\;\cong\;X(q_1,a_1)\;
\oplus\;X(q_3,a_3),$$
then $X(q_2,a_2)\;\cong\;X(q_3,a_3),$ so that $q_2=q_3$
and $a_2=a_3.$
\end{theorem}

Since it is evident that $X(1,0)$ is isomorphic to $\ca$ viewed
as a right $\ca$-module, it follows from the above theorem that
$X(q,0)$ is isomorphic
to the free
$\ca$-module $\ca^q$.
Our previous paper \cite{pacrief} can be used to prove in
a different manner that
$X(q,0)\;\cong\; \ca^q$, and
as mentioned there, filter
functions for dilation of $\br$ by $q$ can be used
to construct module bases for $X(q,0).$
We remark that, as discussed in \cite{Rie1},
$K_0(\ca)\;\cong\;\mathbb Z^2\;=\;\{(q,a):\;q,a\;\in\;\mathbb
Z\}$, consistent with
the above parametrization, so that the positive
cone of
$K_0(\ca)$ is equal to $\mathbb
N\;\times\mathbb
Z\;\cup\;\{(0,0)\}\;\subseteq\;\mathbb Z^2.$ Also,
cancellation fails for complex vector bundles
over tori of dimension
$\geq 5$, and this phenomenon played a crucial role in
our discussion of the (non)-existence of wavelet filter
functions in \cite{pacrief}.

Now Theorem \ref{thmrief} shows that
there is a Hilbert $\ca$-module isomorphism
between $X(q,a)$ and $X(1,a)\;\oplus\; \ca^{q-1}$.
{}From the above observations we then obtain:
\begin{proposition}
\label{thm lpfm}
Let $q\;\in\;\mathbb N,$ and $a\;\in\;\mathbb
Z\;\backslash\{0\}.$
Then there is a module frame in $X(q,a)$ which consists
of $q+1$
elements.
\end{proposition}

It is possible to explicitly write down this module frame
in terms of wavelet filter functions. We leave this
as an exercise for the interested reader.

\section{The ``scaling function"}

In this section we restrict attention to the case of dimension 2,
and to diagonal dilation
matrices. For this case we show how to produce
many non-free projective
$\ca$-submodules $V$
of $\Xi$ which produce projective
multiresolution analyses as defined in Definition \ref{defpr}.
Throughout this section we
assume that $\ca = C(\bt^2)$,
and we fix our dilation matrix $A$ to be of the
form
$A\;=\;\left(\begin{array}{rr}
d_1&0\\
0&d_2
\end{array}\right)$
for $d_i \in \bz$ with $|d_i|>1$ for $i = 1,2$. Thus our earlier
$\d$ will be $\d = (|d_1d_2|)^{-1/2}$. Also, $\Xi$ will
be a subspace of $L^2(\br^2)$.

In the next theorem the function $\s$ acts much like a
traditional scaling function, with condition 2 containing
the analog of the traditional scaling equation.

\begin{theorem}
\label{thmsd}
Fix integers $q\,, a\;\in\;\mathbb Z$ with $q>0$.
Let $\s \in \Xi$ satisfy the further conditions that:
\begin{enumerate}
\item
$
\langle \s,\;\s \rangle_{C(\mathbb
R^2/(\mathbb Z \times q \mathbb Z))}=1,
$
\quad where we use the notation analogous to that of Section 1.
\item There is an $\tilde m \in X(q,(1-d_1d_2)a)$ such that
\[
\s(Ax) = \tilde m(x)\s(x)
\]
for all $x \in \br^2$.
\end{enumerate}
Define ${\mathcal R}:\;X(q,a)\;\rightarrow\;\Xi$ by
the pointwise product
$$
\mcr(F) = \s F.
$$
Then $\mcr$ is an $\ca$-module monomorphism which
has the further properties that:
\begin{enumerate}
\item $\langle {\mathcal R}(F), {\mathcal
R}(G) \rangle_{\ca}\;=\;\langle F,G
\rangle_{\ca},$
\quad for all $F,\;G\;\in\;X(q,a),$
\item ${\mathcal R}(X(q,a))\;\subseteq\;
D{\mathcal R}(X(q,a))$
\end{enumerate}
In particular, ${\mathcal R}(X(q,a))$ is a projective
$\ca$-submodule of $\Xi$.
\end{theorem}

\begin{proof}
A straight-forward calculation, using coset representatives
for $q\bz$ in $\bz$, verifies property (1) above. In particular,
this shows that $\mcr(F)$ is indeed in $\Xi$. Since the module
action is pointwise multiplication on both sides, it is
clear that $\mcr$ is an $\ca$-module homomorphism. From
property (1) it is also clear that $\mcr$ is a monomorphism.

To verify property (2) we must show that if $F \in X(q,a)$ then
$D^{-1}\mcr F \in \mcr X(q,a)$. Recall that
$(D\xi)(x) = \d\xi((A^t)^{-1} x  )$. For our diagonal matrix $A$
we have $A^t = A$, and so $(D^{-1}\xi)(x) = \d^{-1}\xi(Ax)$.
Thus we must show that $(\s F)\circ A = \s G$ for some
$G \in X(q,a)$. Now
\[
((\s F)\circ A)(x) = \s(Ax)F(Ax) = \s(x)\tilde m(x)F(Ax),
\]
and so we define $G$ by $G(x) =  \tilde m(s)F(Ax)$.
Then
\begin{eqnarray*}
G(s+1, t+q) &=& \tilde m(s+1, t+q)F(d_1(s+1), d_2(t+q))  \\
&=& e((1-d_1d_2)as) \tilde m(s,t)e(d_2ad_1s)F(d_1s, d_2t)
= e(as)G(s, t).
\end{eqnarray*}
Thus $G \in X(q,a)$ as desired.
\end{proof}

Our task now is to show the existence of functions $\s$ and
$\tilde m$ satisfying the conditions of the above theorem. These
conditions impose further conditions on $\tilde m$.

\begin{proposition}
Let $\s$ and $\tilde m$ satisfy the conditions of the above
theorem. Let $C$ be a set of coset representatives for
$A(\bz\x q\bz)$ in $\bz\x q\bz$. Then
\[
\sum_{c \in C} |\tilde m(x-A^{-1}c))|^2 = 1.
\]
In other words,
$\<\tilde m, \tilde m\>_{C(\br^2/A^{-1}(\bz\x q\bz))}=1$.
\end{proposition}

\begin{proof}
Using conditions (1) and (2) of the above theorem, we calculate,
for any $x \in \ \br^2$:
\begin{eqnarray*}
1 &=& \<\s,\;\s\>_{C(\bz\x q\bz)}(Ax) =
\sum _{p \in \bz\x q\bz} |\s(Ax-p)|^2
= \sum_{c\in C}\sum_{p \in \bz\x q\bz} |\s(Ax-c -Ap|^2  \\
&=& \sum_c\sum_p |\s(A(x -A^{-1}c -p))|^2
= \sum_c\sum_p|\tilde m(x -A^{-1}c -p)|^2|\s(x -A^{-1}c -p)|^2  \\
&=& \sum_c |\tilde m(x-A^{-1}c)|^2
\<\s, \; \s\>_{C(\bz\x q\bz)}(x-a^{-1}c)
= \sum_c |\tilde m(x-A^{-1}c)|^2
\end{eqnarray*}

\end{proof}

If as coset representatives we choose the integer pairs $(j,k)$
such that $0 \leq j \leq |d_1|-1$ and $0 \leq k \leq |d_2|-1$,
then the above equation can be written as
$$
\sum_{0 \leq j \leq |d_1|-1} \;
\sum_{0 \leq k \leq |d_2|-1}|\tilde{m}(s+\frac{j}{d_1},
t+\frac{kq}{d_2})|^2=1
$$
for all $s,t,\,\in\,\mathbb R$.
When $q=1,$ this equation is closely related to one of
the standard equations that a low-pass filter must
satisfy in ordinary wavelet theory corresponding to
dilation by the matrix $A$, except that now
$\tilde{m}\;\in\;X(q,(1-d_1d_2)a).$

Following the pattern for ordinary scaling functions, we seek
to define $\s$ by
$$
\s(s,t)\;=\;\prod_{i=1}^{\infty}\tilde{m}(A^{-i}(s,t)).
$$
This can be done directly, but we choose to use facts from
ordinary wavelet theory to by-pass having to deal directly with
this infinite product. But our discussion below is
equivalent to finding conditions on $\tilde m$ which assure that
this product converges to a function having the properties
which we need.

As usual, we will insist that
$\tilde m(0, 0) = 1$. Then, as in ordinary wavelet theory, from the
equation in the above proposition it follows
that $\tilde m(0, kq/d_2) = 0$ for $1 \leq k \leq |d_2|-1$.
We will seek $\tilde m$ of the form
$\tilde{m}(s,t)\;=\;J_{q,(1-d_1d_2)a,\b}\;m$ where
$m$ is an ordinary low-pass filter function for dilation by $A$,
but for translation by $\bz \x q\bz$ instead of $\bz^2$
(so that $m \in X(q, 0)$).  But even more,
we take $m$
to be a tensor product of ordinary one-dimensional low-pass
filters, $m_1$ and $m_2$, where $m_1$ will be for dilation
by $d_1$, while $m_2$ will be for dilation by $d_2$ but also
for translation by $q\bz$ instead of $\bz$. (Again, we
normalize
our low-pass filters to have value $1$ at $0.$)
According to
theorem 3.2 of \cite{BS}, and especially its claim 3.3, we
can choose $m_1$ and $m_2$ to be smooth, and to give
corresponding scaling functions $\var_1$ and $\var_2$
which are Schwartz functions.

We will define $\var$ by $\var(s,t) =
\var_1(s)\var_2(t)$ for $s, \, t \in \br$. Then $\var$ itself
is a Schwartz function, and so is in $\Xi$ according to
Corollary \ref{cor1s}.
Clearly $\var$ is the ordinary scaling function for $m$, and
$\langle \var,\;\var \rangle_{C(\mathbb
R^2/(\mathbb Z \times q \mathbb Z))}=1$.
We will have $m_2(0) = 1$,
so that $m_2(q/d_2) = 0$ as usual, and thus we will have
$m(s, q/d_2) = 0$ for all $s$. Thus we will take
$\b = q/d_2$ so that
$\tilde{m}\;=\;J_{q,(1-d_1d_2)a, q/d_2}m$ will be a continuous
function.

  From the infinite product indicated above for $\s$ it is
easy to guess that $\s$ will be of the form $\s = J\var$
for a suitable function $J$. We can expect that $J$ will not
be continuous, but it should be bounded and
have its discontinuities where
$\var$ takes value $0$.
Now as usual $\var$ satisfies the scaling equation
\[
\var(Ax) = m(x)\var(x),
\]
that is,
\[
\var(d_1s, d_2t) = m(s,t)\var(s,t).
\]
Since $m(s, (q/d_2) + nq) = 0$ for any $n \in \bz$, we find,
upon iterating,
that $\var$ will take value $0$ on the vertical lines with
$t = (1+d_2n)d_2^{j-1}q$ for $j \geq 1$ and $n \in \bz$.
Notice that all of these values of $t$ are integers not equal
to $0$.
For simplicity of exposition we will ignore in the next
proposition what happens on these vertical lines.

\begin{proposition}
Define on $\br^2$ a function, $J$, by
\[ J(x) =\prod_{j=1}^{\infty} J_{q,(1-d_1d_2)a, q/d_2}(A^{-j}(x)).
\]
This product converges uniformly on any compact subset of
$\br^2$ which does not meet the vertical lines
$(s, (1+d_2n)d_2^{j-1}q)$ . Thus it
is a continuous function except possibly on those vertical
lines. Furthermore, $|J(x)| = 1$ except possibly on those
vertical lines. Finally,
\[
J(Ax) = J_{q,(1-d_1d_2)a, q/d_2}(x)J(x)
\]
for $x$ not on those vertical lines.
\end{proposition}

\begin{proof}
For any $j \geq 1$ we have\[
J_{q,(1-d_1d_2)a, q/d_2}(A^{-j}(s,t)) =
J_{q,(1-d_1d_2)a, q/d_2}(d_1^{-j}s, d_2^{-j}t)
\]
for all $x = (s,t)$. So the discontinuities of this function
occur exactly where $d_2^{-j}t = (q/d_2) + nq$, which gives
exactly the vertical lines found above. The function
$J_{q,(1-d_1d_2)a, q/d_2}$ is infinitely differentiable at
$(0,0)$ and has value $1$ there,
and so the usual proof for the construction of
scaling functions works here too, when applied to points not
on the vertical lines. See, for example, section 2.3
of \cite{HW}. A simple standard argument then verifies the
above scaling equation for $J$. Further examination using the
fact that $|J_{q,(1-d_1d_2)a, q/d_2}| = 1$ shows that $|J|=1$.
\end{proof}

\begin{theorem}
With notation as above, set $\s = J\var$. Then $\s$ satisfies
the conditions of Theorem \ref{thmsd}.
\end{theorem}

\begin{proof}
Because $\var$ is continuous, and because $J$ is bounded and
its discontinuities
occur exactly where $\var$ takes value $0$,
it follows that $\s$ is continuous.
Furthermore, for any $x \in \br^2$ we have
\[
\s(Ax) = J(Ax)\var(Ax) = J_{q,(1-d_1d_2)a, q/d_2}(x)J(x)m(x)\var(x)
= \tilde m(x)\s(x).
\]
Finally, since $|J_{q,(1-d_1d_2)a, q/d_2}| = 1$, we have
$\<\s, \s\>_\ca = \<\var, \var\>_\ca$ so that $\s \in \Xi$.
\end{proof}

\begin{theorem}
\label{thm4diag}
For our given dilation matrix $A$ and
for a given $X(q,a)$ there are projective multiresolution analyses
$\{V_j\}$ such that $V_0 \cong X(q,a)$.
\end{theorem}

\begin{proof}
We have obtained a function $\s$ satisfying
the conditions of Theorem \ref{thmsd}, from which we obtain
the corresponding module monomorphism $\mcr$ of $X(q,a)$
into $\Xi$. We set $V= V_0 = \mcr (X(q,a)$. From Theorem \ref{thmsd}
we see that $V \subseteq D(V)$. For each $j \in \bz$
we set $V_j = D(V)$. Then from Propositions \ref{prop2i}
and \ref{prop2d}
we see that the family $\{V_j\}$
of subspaces of $\Xi$ is a projective multiresolution
analysis as defined in Definition \ref{defpr}.
\end{proof}

\section{The structure of the wavelet module}

{}From the construction of the projective multiresolution analyses
just given we can then construct the corresponding wavelet
modules $W_j$ as discussed before Theorem \ref{thm2w}. We now
come another of our main results, which is the identification
of the isomorphism class of $W_0$ as projective $\ca$-module.
In order to
do this we must first identify the isomorphism class of $V_1$
as $\ca$-module.

Recall that for our given ``scaling function'' $\s$ we have
\[
V_1 = \{(\s F)\circ B: F \in X(q,a)\}.
\]
Set $Y_A(q, a) = \{F\circ B: \; F \in X(q,a)\}$, with its
evident structure as an $\ca$-module by pointwise multiplication.
The mapping $(F\circ B) \mapsto (\s\circ B)(F\circ B)$ is clearly
an $\ca$-module isomorphism from $Y_A(q,a)$ onto $V_1$. (So
$Y_A(q,a)$ must be projective.) Thus it
suffices for us to determine the isomorphism class of $Y_A(q,a)$.

Now $(F\circ B)(s,t) = F(s/d_1, t/d_2)$, and simple calculations
show that equivalently we have
\[
Y_A(q,a) = \{G \in
C_b(\br^2): G(s-d_1, t) = G(s,t) \quad\mathrm{and}\quad G(s,t-d_2q)
= e((a/d_1)s)G(s,t)\}.
\]
Notice that $a/d_1$ need not be an integer. To deal with this we
let $c$ be the greatest common divisor of $a$ and $d_1$, chosen
so that $c>0$, and we define $\hat a$ and $\hat d_1$ by
$a=\hat a c$ and $d_1 = \hat d_1 c$.

For any $p, q, a \in \bz$ with $p\neq 0$ and $q\neq 0$ set
\[
Z(p;q,a) = \{F \in C_b(\br^2): F(s-p, t) = F(s,t)
\quad\mathrm{and}\quad F(s, t-q) = e(as)F(s,t)\}.
\]

\begin{lemma}
\label{lem5y}
The $\ca$-module $Z(d_1;\hat d_1 d_2q, \hat a)$ is isomorphic
to the direct sum of $|\hat d_1|$ copies of $Y_A(q,a)$, which
we denote by \ $|\hat d_1| Y_A(q,a)$.
\end{lemma}

\begin{proof}
Define an automorphism, $\a$, of $Z(d_1;\hat d_1 d_2q, \hat a)$
by
\[
(\a F)(s,t) = \bar e((\hat a/\hat d_1)s)F(s, t-d_2q).
\]
It is easily checked that this is indeed a module automorphism,
and that it is of order $|\hat d_1|$. That is, the cyclic
group of order $|\hat d_1|$ acts on the module. Then
$Z(d_1;\hat d_1 d_2q, \hat a)$ decomposes into
the direct sum of its ``isotypic components'' for this action.
The fixed-point
submodule clearly consists of the $F$'s in
$Z(d_1;\hat d_1 d_2q, \hat a)$ satisfying
\[
F(s, t-d_2q) = e((a/d_1)s)F(s,t)\}.
\]
That is, the fixed-point submodule actually is $Y_A(q,a)$.
For any $k\in \bz$ with $0 \leq k \leq |\hat d_1| - 1$ the
$k$-th isotypic component consists of the $F$'s which
satisfy $\a(F) = e(k/\hat d_1)F$. But for such an $F$ set
$F_0(s, t) = e(kt/(\hat d_1d_2q))F(s,t)$. Then it is easily checked
that $F_0 \in Y_A(q,a)$, and that this gives an $\ca$-module
isomorphism of the $k$-th isotypic component onto $Y_A(q,a)$.
\end{proof}

So we see that we must now determine the isomorphism class
of modules of the form $Z(p;q,a)$.

\begin{proposition}
\label{prop5z}
For and $p, q, a \in \bz$ with $p \neq 0$ and $q \neq 0$ the
$\ca$-module $Z(p;q,a)$ is isomorphic to the direct sum of
$|p|$ copies of $X(q,a)$, and so to $X(pq, pa)$.
\end{proposition}

\begin{proof}
Define an automorphism, $\a$, of $Z(p;q,a)$ by
$(\a F)(s,t) = F(s-1, t)$. Then $\a$ has order $|p|$, so that
the cyclic group of order $|p|$ acts, and $Z(p;q,a)$ decomposes
into the direct sum of its isotypic components. The fixed-point
submodule is clearly $X(q,a)$. For any $k\in \bz$ with
$0 \leq k \leq  |p|- 1$ the
$k$-th isotypic component consists of the $F$'s which
satisfy $\a(F) = e(k/p)F$. But for such an $F$ set
$F_0(s, t) = \bar e(ks/p)F(s,t)$. It is easily checked that
$F_0 \in Y(q,a)$, and that this gives a module isomorphism
of the $k$-th isotypic component onto $Y(q,a)$. Thus
$Z(p;q,a)$ is isomorphic to $|p|$ copies of $Y(q,a)$, which
we can denote by $|p|Y(q,a)$. But by
Theorem \ref{thmrief} this in turn is isomorphic to $X(pq, pa)$.
\end{proof}

\begin{theorem}
With notation as above,  as $\ca$-modules we have
\[V_1 \cong Y_A(q,a)\cong X((\det(A))q, a).
\]
\end{theorem}

\begin{proof}
We saw above that $V_1 \cong Y_A(q,a)$.
Upon applying successively Proposition \ref{lem5y},
Proposition \ref{prop5z}, and Theorem \ref{thmrief}, we
find that
\[
|\hat d_1|Y_A(q,a) \cong Z(d_1;\hat d_1d_2q, \hat a)
\cong |d_1|X(\hat d_1d_2q, \hat a)
\cong |\hat d_1|X(d_1d_2q, c\hat a).
\]
We know that $Y_A(q,a)$ is a projective module, and so
by the cancellation property stated as part of Theorem
\ref{thmrief} we obtain the desired conclusion.
\end{proof}

\begin{theorem}
\label{thm5wavmod}
Let $W_0$ be the wavelet space for the projective multiresolution
analysis based as above on the projective  module $V$ isomorphic
to $X(q,a)$, for the diagonal dilation matrix $A$. Then
\[
W_0 \cong X\big((|\det(A)|- 1)q,(\mathrm{sign}(\det(A))-1)a\big).
\]
Thus if $\det(A) > 0$ then $W_0$ is a free module, while if
$\det(A) < 0$ then $W_0$ is not a free module as long as
$a \neq 0$.
\end{theorem}

\begin{proof}
  From the comment before Theorem \ref{thmrief} we have
\[
X((\det(A))q, a) = X(|\det(A)|q, \mathrm{sign}(\det(A))a).
\]
  From the above theorem we see that
\[
W_0 \oplus X(q,a) \cong
X\big(|\det(A)|q,\mathrm{sign}(\det(A))a\big).
\]
The desired conclusion then follows from Theorem \ref{thmrief}.
\end{proof}

Both Theorem \ref{thm4diag} and Theorem \ref{thm5wavmod} deal with
the case where the
dilation matrix is a $2\times 2$ diagonal matrix, and it is natural to ask
whether or not it is possible to
build a projective multiresolution analysis whose initial
module $V$ is isomorphic to any $C(\mathbb T^2)$-module
for any $2\times 2$ integer dilation matrix.
As mentioned before, Bownik and Speegle have shown that for
any such dilation matrix, there exist scaling functions
whose Fourier transforms are continuous and compactly
supported \cite{BS}. For this reason, we conjecture that
such a construction is possible.  However, to date, we have
only been able to perform this construction for matrices that
are similar to diagonal matrices, and for the dilation
matrix $A\;=\;\left(\begin{array}{rr}
1&-1\\
1&1
\end{array}\right)$ corresponding to the quincunx lattice.


\end{document}